\documentclass[12pt]{amsart}
\usepackage{mathrsfs}
\usepackage{}
\usepackage{amsmath}
\usepackage{amsfonts}
\usepackage{amssymb}
\usepackage[all,cmtip]{xy}           
\usepackage{shuffle}
\usepackage{caption}
\usepackage{bbding}
\usepackage{txfonts}
\usepackage[shortlabels]{enumitem}
\usepackage{ifpdf}
\ifpdf
\usepackage[colorlinks,final,backref=page,hyperindex]{hyperref}
\else
\usepackage[colorlinks,final,backref=page,hyperindex,hypertex]{hyperref}
\fi
\usepackage{tikz}
\usepackage[active]{srcltx}
\usepackage{float}

\makeatletter

\newtheorem{theorem}{Theorem}[section]

\newtheorem{lemma}[theorem]{Lemma}

\newtheorem{prop-def}{Proposition-Definition}[section]

\theoremstyle{definition}
\newtheorem{defn}[theorem]{Definition}

\newtheorem{remark}[theorem]{Remark}

\newcommand{\nc}{\newcommand}

\newcommand {\emptycomment}[1]{}

\nc{\delete}[1]{{}}
\nc{\mmargin}[1]{}

\nc{\mlabel}[1]{\label{#1}}  
\nc{\mcite}[1]{\cite{#1}}  
\nc{\mref}[1]{\ref{#1}}  
\nc{\meqref}[1]{\eqref{#1}}  
\nc{\mbibitem}[1]{\bibitem{#1}} 

\newcommand{\bk}{{\mathbf{k}}}

\nc{\vep}{\varepsilon}
\nc{\bin}[2]{ (_{\stackrel{\scs{#1}}{\scs{#2}}})}  
\nc{\binc}[2]{(\!\! \begin{array}{c} \scs{#1}\\
		\scs{#2} \end{array}\!\!)}  
\nc{\bincc}[2]{  ( {\scs{#1} \atop
		\vspace{-1cm}\scs{#2}} )}  
\nc{\oline}[1]{\overline{#1}}
\nc{\mapm}[1]{\lfloor\!|{#1}|\!\rfloor}
\nc{\bs}{\bar{S}}
\nc{\cast}{{\,\mbox{\raisebox{.8pt}{$\scriptstyle \circledast$}}\,}}
\nc{\la}{\longrightarrow}
\nc{\ot}{\otimes}
\nc{\rar}{\rightarrow}
\nc{\dar}{\downarrow}
\nc{\dap}[1]{\downarrow \rlap{$\scriptstyle{#1}$}}
\nc{\defeq}{\stackrel{\rm def}{=}}
\nc{\dis}[1]{\displaystyle{#1}}
\nc{\dotcup}{\ \displaystyle{\bigcup^\bullet}\ }
\nc{\hcm}{\ \hat{,}\ }
\nc{\hts}{\hat{\otimes}}
\nc{\hcirc}{\hat{\circ}}
\nc{\lleft}{[}
\nc{\lright}{]}
\nc{\curlyl}{\left \{ \begin{array}{c} {} \\ {} \end{array}
	\right .  \!\!\!\!\!\!\!}
\nc{\curlyr}{ \!\!\!\!\!\!\!
	\left . \begin{array}{c} {} \\ {} \end{array}
	\right \} }
\nc{\longmid}{\left | \begin{array}{c} {} \\ {} \end{array}
	\right . \!\!\!\!\!\!\!}
\nc{\ora}[1]{\stackrel{#1}{\rar}}
\nc{\ola}[1]{\stackrel{#1}{\la}}
\nc{\scs}[1]{\scriptstyle{#1}} \nc{\mrm}[1]{{\rm #1}}
\nc{\dirlim}{\displaystyle{\lim_{\longrightarrow}}\,}
\nc{\invlim}{\displaystyle{\lim_{\longleftarrow}}\,}
\nc{\dislim}[1]{\displaystyle{\lim_{#1}}} \nc{\colim}{\mrm{colim}}
\nc{\mvp}{\vspace{0.3cm}} \nc{\tk}{^{(k)}} \nc{\tp}{^\prime}
\nc{\ttp}{^{\prime\prime}} \nc{\svp}{\vspace{2cm}}
\nc{\vp}{\vspace{8cm}}
\nc{\modg}[1]{\!<\!\!{#1}\!\!>}
\nc{\intg}[1]{F_C(#1)}
\nc{\lmodg}{\!<\!\!}
\nc{\rmodg}{\!\!>\!}
\nc{\cpi}{\widehat{\Pi}}
\nc{\labs}{\mid\!}
\nc{\rabs}{\!\mid}
\nc{\btr}{\blacktriangleright}

\nc{\ad}{\mrm{ad}}
\nc{\rRB}{\mathsf{rRB}}
\nc{\cocrRB}{\mathsf{cocrRB}}
\nc{\PH}{\mathsf{PH}}
\nc{\cocPH}{\mathsf{cocPH}}
\nc{\ann}{\mrm{ann}}
\nc{\Ad}{\mrm{Ad}}
\nc{\Aut}{\mrm{Aut}}
\nc{\Der}{\mrm{Der}}
\nc{\Sym}{\mrm{Sym}}
\nc{\br}{\mrm{bre}}
\nc{\can}{\mrm{can}}
\nc{\Cont}{\mrm{Cont}}
\nc{\rchar}{\mrm{char}}
\nc{\cok}{\mrm{coker}}
\nc{\de}{\mrm{dep}}
\nc{\dtf}{{R-{\rm tf}}}
\nc{\dtor}{{R-{\rm tor}}}

\nc{\Dif}{\mrm{Diff}}
\nc{\Div}{\mrm{Div}}
\nc{\End}{\mrm{End}}
\nc{\Ext}{\mrm{Ext}}
\nc{\Fil}{\mrm{Fil}}
\nc{\Fr}{\mrm{Fr}}
\nc{\Frob}{\mrm{Frob}}
\nc{\Gal}{\mrm{Gal}}
\nc{\GL}{\mrm{GL}}
\nc{\Gr}{\mrm{Gr}}
\nc{\Hom}{\mrm{Hom}}
\nc{\Hoch}{\mrm{Hoch}}
\nc{\hsr}{\mrm{H}}
\nc{\hpol}{\mrm{HP}}
\nc{\id}{\mrm{id}}
\nc{\im}{\mrm{im}}
\nc{\inv}{\mrm{inv}}
\nc{\Id}{\mrm{Id}}
\nc{\ID}{\mrm{ID}}
\nc{\Irr}{\mrm{Irr}}
\nc{\incl}{\mrm{incl}}
\nc{\length}{\mrm{length}}
\nc{\NLSW}{\mrm{NLSW}}
\nc{\Lie}{\mrm{Lie}}
\nc{\mchar}{\rm char}
\nc{\mpart}{\mrm{part}}
\nc{\ql}{{\QQ_\ell}}
\nc{\qp}{{\QQ_p}}
\nc{\rank}{\mrm{rank}}
\nc{\rcot}{\mrm{cot}}
\nc{\rdef}{\mrm{def}}
\nc{\rdiv}{{\rm div}}
\nc{\rtf}{{\rm tf}}
\nc{\rtor}{{\rm tor}}
\nc{\res}{\mrm{res}}
\nc{\SL}{\mrm{SL}}
\nc{\Spec}{\mrm{Spec}}
\nc{\tor}{\mrm{tor}}
\nc{\Tr}{\mrm{Tr}}
\nc{\tr}{\mrm{tr}}
\nc{\wt}{\mrm{wt}}

\nc{\bfk}{{\bf k}}
\nc{\bfone}{{\bf 1}}
\nc{\bfzero}{{\bf 0}}
\nc{\detail}{\marginpar{\bf More detail}
	\noindent{\bf Need more detail!}
	\svp}
\nc{\gap}{\marginpar{\bf Incomplete}\noindent{\bf Incomplete!!}
	\svp}
\nc{\FMod}{\mathbf{FMod}}
\nc{\Int}{\mathbf{Int}}
\nc{\Mon}{\mathbf{Mon}}
\nc{\remarks}{\noindent{\bf Remarks: }}
\nc{\Rep}{\mathbf{Rep}}
\nc{\Rings}{\mathbf{Rings}}
\nc{\Sets}{\mathbf{Sets}}
\nc{\Diff}{\mathbf{Diff}}
\nc{\Inte}{\mathbf{Inte}}
\nc{\U}{\mathrm{U}}

\nc{\BA}{{\mathbb A}}   \nc{\CC}{{\mathbb C}}
\nc{\DD}{{\mathbb D}}   \nc{\EE}{{\mathbb E}}
\nc{\FF}{{\mathbb F}}   \nc{\GG}{{\mathbb G}}
\nc{\HH}{{\mathbb H}}   \nc{\LL}{{\mathbb L}}
\nc{\NN}{{\mathbb N}}   \nc{\PP}{{\mathbb P}}
\nc{\QQ}{{\mathbb Q}}   \nc{\RR}{{\mathbb R}}
\nc{\TT}{{\mathbb T}}   \nc{\VV}{{\mathbb V}}
\nc{\ZZ}{{\mathbb Z}}   \nc{\TP}{\widetilde{P}}

\nc{\cala}{{\mathcal A}}    \nc{\calc}{{\mathcal C}}
\nc{\cald}{\mathcal{D}}     \nc{\cale}{{\mathcal E}}
\nc{\calf}{{\mathcal F}}    \nc{\calg}{{\mathcal G}}
\nc{\calh}{{\mathcal H}}    \nc{\cali}{{\mathcal I}}
\nc{\call}{{\mathcal L}}    \nc{\calm}{{\mathcal M}}
\nc{\caln}{{\mathcal N}}    \nc{\calo}{{\mathcal O}}
\nc{\calp}{{\mathcal P}}    \nc{\calr}{{\mathcal R}}

\nc{\cals}{{\mathcal S}}    \nc{\calt}{{\Omega}}
\nc{\calv}{{\mathcal V}}    \nc{\calw}{{\mathcal W}}
\nc{\calx}{{\mathcal X}}

\nc{\fraka}{{\mathfrak a}}
\nc{\frakb}{\mathfrak{b}}
\nc{\frakg}{{\frak g}}
\nc{\frakl}{{\frak l}}
\nc{\fraks}{{\frak s}}
\nc{\frakB}{{\frak B}}
\nc{\frakm}{{\frak m}}
\nc{\frakM}{{\frak M}}
\nc{\frakp}{{\frak p}}
\nc{\frakW}{{\frak W}}
\nc{\frakX}{{\frak X}}
\nc{\frakS}{{\frak S}}
\nc{\frakA}{{\frak A}}
\nc{\frakx}{{\frakx}}

\nc{\ynr}[1]{\textcolor{orange}{\underline{Yunnan:}#1 }}

\nc{\lir}[1]{\textcolor{red}{\underline{Li:}#1 }}

\delete{
	\input cyracc.def

	\newtheorem{theorem}{Theorem}[section]

	\theoremstyle{definition}

	\theoremstyle{remark}
	\newtheorem{remark}[theorem]{Remark}

	\allowdisplaybreaks
	
	\numberwithin{equation}{section}
	
}

\begin{document}

\title[Matched pairs of actions on the Kac-Paljutkin algebra $H_8$]{Matched pairs of actions on the Kac-Paljutkin algebra $H_8$}

\author{Yongyue Xiao}
\address{School of Mathematics and Information Science, Guangzhou University,
Guangzhou 510006, China}
\email{2112315047@e.gzhu.edu.cn}

\author{Yunnan Li}
\address{School of Mathematics and Information Science, Guangzhou University,
Guangzhou 510006, China}
\email{ynli@gzhu.edu.cn}

\begin{abstract}
The notion of matched pair of actions on a Hopf algebra generalizes the braided group construction of Lu, Yan and Zhu, and efficiently provides Yang-Baxter operators.
In this paper, we classify matched pairs of actions on the Kac-Paljutkin Hopf algebra $H_8$. Through calculations, we obtain 6 matched pairs of actions on $H_8$. Based on such a classification result, we find that four of them can be derived from the  coquasitriangular structures of $H_8$, while the other two can not. Furthermore, we discover that the Yang-Baxter operators associated to exactly these two distinguished matched pairs of actions are involutive.
\end{abstract}

\keywords{matched pair, the Kac-Paljutkin Hopf algebra, Yang-Baxter operator
\\
\qquad 2020 Mathematics Subject Classification. 16T05, 16T25}

\maketitle

\tableofcontents

\allowdisplaybreaks

\section{Introduction}
The concept of a matched pair of Hopf algebras was first introduced by Singer \cite{Si} in the graded case in 1972. This idea was later extended to the ungraded case in \cite{Ta}, and Takeuchi also introduced the notion of a matched pair of groups.
Subsequently, Majid presented the modern definition of a matched pair of Hopf algebras \cite{Ma2,Ma} with the aim of providing an elegant description for the Drinfeld double of a finite dimensional Hopf algebra. As detailed  in \cite{Ma}, the definition of a matched pair of Hopf algebras is as follows.
\begin{defn}[{\cite[\S~7.2]{Ma}}]
A {\bf matched pair of Hopf algebras} is a quadruple $(H,K,\rightharpoonup,\leftharpoonup)$, where $H$ and $K$
are Hopf algebras, $\rightharpoonup$ is a left $H$-module coalgebra action on $K$, $\leftharpoonup$ is a right $K$-module coalgebra action on $H$ such that\vspace{-.5em}
\begin{eqnarray}
\label{eq:MP1}
x\rightharpoonup ab&=&(x_1\rightharpoonup a_1)((x_2\leftharpoonup a_2)\rightharpoonup b),\\
\label{eq:MP2}
x\rightharpoonup 1_K&=&\varepsilon_H(x)1_K,\\
\label{eq:MP3}
xy\leftharpoonup a&=&(x\leftharpoonup(y_1\rightharpoonup a_1))(y_2\leftharpoonup a_2),\\
\label{eq:MP4}
1_H\leftharpoonup a&=&\varepsilon_K(a)1_H,\\
\label{eq:MP5}
(x_1\rightharpoonup a_1) \otimes (x_2\leftharpoonup a_2)  &=&
(x_2\rightharpoonup a_2) \otimes (x_1\leftharpoonup a_1)
\end{eqnarray}
for $x,y\in H$ and $a,b\in K$. Here we compress Sweedler's notation of the comultiplication $\Delta$ of a coalgebra $C$ as
$$\Delta(x)=x_1\otimes x_2$$
for simplicity. For other basic notions of Hopf algebras, we refer to the textbooks~\mcite{Mon}.
\end{defn}

Following this, Agore, Bontea and Militar applied matched pairs of Hopf algebras to solve factorization problem for Hopf algebras in \cite{ABM}, where they also classified all matched pairs $(H_4,\bk[C_n],\rightharpoonup,\leftharpoonup)$ and $(H_4,\bk[C_2\times C_2],\rightharpoonup,\leftharpoonup)$.
Later, Agore classified all matched pairs of two Taft algebras~\cite{Ag}.
The classification of matched pairs for various other Hopf algebras of low dimension can be found in works such as \cite{Bo,LNW}.

Recently, Ferri and Sciandra further introduced the notion of matched pairs of actions on Hopf algebras, as a certain subclass of matched pairs of Hopf algebras with compatible actions.
\begin{defn}[\cite{FS}]
Given a Hopf algebra $H=(H,\cdot\,,1,\Delta,\vep,S)$, a matched pair $(H,H,\rightharpoonup,\leftharpoonup)$ is called a {\bf matched pair of actions} on $H$, when
\begin{eqnarray}
\label{eq:MP*}
xy&=&(x_1\rightharpoonup y_1)(x_2\leftharpoonup y_2),\quad\forall x,y\in H.
\end{eqnarray}
Equivalently, we have
\begin{eqnarray}\label{eq:l-action}
x\rightharpoonup y&=&x_1y_1S(x_2\leftharpoonup y_2),\quad\forall x,y\in H.
\end{eqnarray}
\begin{eqnarray}\label{eq:r-action}
x\leftharpoonup y&=&S(x_1\rightharpoonup y_1)x_2y_2,\quad\forall x,y\in H.
\end{eqnarray}
A matched pair of actions on $H$ will be abbreviated as $(H,\rightharpoonup,\leftharpoonup)$.
\end{defn}
The work of Ferri and Sciandra generalizes the seminal construction of braided group by Lu, Yan and Zhu in~\cite{LYZ}, and also the counterpart for cocommutative Hopf algebras originally studied by Angiono, Galindo, and Vendramin in \cite{AGV}. It is worth noting that in both \cite{LYZ} and \cite{AGV}, matched pairs of groups or cocommutative Hopf algebras were used to construct solutions to the Yang-Baxter equation.
In the recently updated work \cite{GGV}, Guccione, Guccione and Valqui further pointed out that matched pairs of actions on arbitrary Hopf algebras equivalently produce braiding operators satisfying the braid relation. So in our opinion, it is meaningful to study matched pairs of actions on Hopf algebras for finding Yang-Baxter operators, and
certain equivalent conditions have been provided in the preceding paper~\cite{Li} for such kind of Yang-Baxter operators being involutive.
Particularly for certain non-commutative and non-cocommutative Hopf algebras of low dimension, we aim to classify their matched pairs of actions in order to obtain concrete Yang-Baxter operators.

In \cite{FS}, Ferri and Sciandra also pointed out that coquasitriangular bialgebras produce matched pairs of actions and provided concrete examples to illustrate this result.
For instance, the smallest example $H_4$, the 4-dimensional Sweedler Hopf algebra, has a unique matched pair of actions. In this paper, we further investigate matched pairs of actions on the Kac-Paljutkin algebra $H_8$. It is the unique 8-dimensional non-commutative and non-cocommutative semisimple Hopf algebra arising from exact factorizations of
finite groups \cite{KP,Mas}.
According to \cite{St}, there are 7 non-isomorphic 8-dimensional Hopf algebras, namely $A_{C_2}$, $A_{C_2\times C_2}$, $A'_{C_4}$, $A''_{C_4}$, $A'''_{C_4,q}$, $(A''_{C_4})^*$ and $H_8$, excluding group algebras and their duals. After dealing with the most complicated example $H_8$ among them, one can discuss other non-semisimple cases in a similar manner (see e.g.~\cite[Example~4.3]{Li}).

So far the Kac-Paljutkin algebra $H_8$ has been widely studied from different aspects. For instance, all its coquasitriangular structures were determined by Suzuki in his study of a more general kind of finite-dimensional cosemisimple Hopf algebras $A^{\nu\lambda}_{NL}$~\cite{Su}. All quasitriangular structures on $H_8$ were explicitly described by Wakui in~\cite{Wa}; see also~\cite{Pa,ZL} for a broader class of
semisimple Hopf algebras $H_{2n^2}$ including the Kac-Paljutkin algebra. The twisted Frobenius-Schur indicators for $H_8$ were given by Sage and Vega in \cite{SV}. All simple modules for the Drinfeld double of $H_8$ were classified by Hu and Zhang in \cite{HZ1}.
All finite-dimensional Nichols algebras and then Hopf algebras over $H_8$ were classified by Shi in~\cite{Shi}.

The paper is organized as follows.
In Section \ref{sec:H8}, in order to effectively find all the matched pairs of actions on $H_8$, we first propose a classification strategy in Section \ref{sec:st}. Then, we apply the strategy and obtain 6 matched pairs of actions on $H_8$ (Theorem~\ref{thm:mpa_H_8}) in Section \ref{sec:cl}. Based on such a classification result, we discuss the relationship between matched pairs of actions on $H_8$ and its coquasitriangular structures in Section \ref{sec:ir1}. Our finding reveals that four of the matched pairs of actions on $H_8$ can be derived from its own coquasitriangular structures, while the remaining two can not (Theorem~\ref{thm:mpa_H_8_coqt}). Based on the previous discussion, we also explore the relationship between matched pairs of actions on $H_8$ and Yang-Baxter operators in Section \ref{sec:ir2}. Finally, we discover that exactly the two Yang-Baxter operators associated to the matched pairs of actions on $H_8$, which are not derived from its coquasitriangular structures, are involutive (Theorem~\ref{thm:mpa_H_8_ybo}).

\noindent
{\bf Convention.}
In this paper, we fix an algebraically closed ground field $\bk$ of characteristic 0.
All the objects under discussion, including vector spaces, algebras and tensor products, are taken over $\bk$ by default.

\section{Matched pairs of actions on the Kac-Paljutkin Hopf algebra}\label{sec:H8}

In this section, we classify all matched pairs of actions on the Kac-Paljutkin Hopf algebra $H_8$, the unique non-commutative and non-cocommutative semisimple
Hopf algebra of dimension $8$.

The algebra $H_8$ is generated by 3 elements $g$, $h$ and $z$ with the relations
$$g^{2} =h^{2}= 1,\quad z^{2} = \frac{1}{2}(1+g+h-gh), \quad gh = hg,\quad zg=hz,\quad zh=gz.$$
A linear basis of $H_8$ is given by $\{1, g, h, gh, z, gz,hz, ghz\}$. It is easy to verify that the element $z$ is invertible with $z^4=1$. The coalgebra structure and the antipode of $H_8$ are defined by
$$
\Delta(g) = g \otimes g, \quad \Delta(h) = h \otimes h,
\quad \Delta(z) =J(z\otimes z),$$
$$
\varepsilon(g) = \varepsilon(h) = \varepsilon(z) = 1, \quad S(g) = g, \quad S(h) = h, \quad S(z) = z.
$$
where $J=\frac{1}{2} (1 \otimes 1 + g \otimes 1+ 1 \otimes h- g \otimes h )$.
In particular, the group of group-like elements in $H_8$
$$G(H_8)=\{1,g,h,gh\},$$
isomorphic to the Klein 4-group.

\subsection{Strategy for classification of matched pairs of actions on $H_8$}\label{sec:st}
\
\newline

\vspace{-.7em}

First recall the following characterization of matched pairs of actions on a Hopf algebra given in \cite[Theorem~4.1]{Li}; see also \cite[Theorem~5.27]{GGV}. It will be useful for our later classification.
\begin{theorem}\label{thm:construct-mp}
Let $H$ be a Hopf algebra, and $\rightharpoonup$ be a left $H$-module coalgebra action on itself. Define linear map $\leftharpoonup:H\otimes H\to H$ by Eq.~\eqref{eq:r-action}.
If $\leftharpoonup$ is a right $H$-module coalgebra action,
then $(H,\rightharpoonup,\leftharpoonup)$
is a matched pair of actions on $H$.
\end{theorem}

According to Theorem~\ref{thm:construct-mp},
in order to classify matched pairs of actions on $H_8$, it is enough to calculate the interactions among its 8 basis elements by all Eqs.~\eqref{eq:MP1}--\eqref{eq:MP*} to determine left module coalgebra actions $\rightharpoonup$, and then check that linear maps $\leftharpoonup$ obtained from Eq.~\eqref{eq:r-action} are right module coalgebra actions.

To effectively find all the matched pairs of actions on $H_8$, we propose the following strategy, which is summarized with inspiration drawn from \cite{Bo} and \cite{LNW}.

First, we consider the interactions among group-like elements because the coproduct of them is simple. To make the calculations more efficient, we utilize Eq.~\eqref{eq:r-action} to relate the left and right actions among the group-like elements. In this way, once we obtain the left action, the right action is determined accordingly.

Next, we consider the interactions between the group-like elements and $z,\,gz,\,hz,\,ghz$. Although the coproduct of $z,\,gz,\,hz,\,ghz$ is complicated, fortunately, we have found that the left and right actions of $z,\,gz,\,hz,\,ghz$ on group-like elements still result in group-like elements, which significantly reduces the workload of calculations. Therefore, to complete this step, we only need to use Eqs.~\eqref{eq:MP1},~\eqref{eq:MP3},~\eqref{eq:l-action} and~\eqref{eq:r-action} to calculate the left and right actions of the group-like elements on $z,\,gz,\,hz,\,ghz$.

Lastly, we consider the interactions among $z,\,gz,\,hz,\,ghz$. The implementation of this step is much more difficult than the aforementioned steps. But in fact, we realize that to determine the interactions among these elements, we only need to calculate the key elements $z\rightharpoonup z,\,z\rightharpoonup gz,\,z\rightharpoonup hz,\,z\rightharpoonup ghz,\,z\leftharpoonup z,\,gz\leftharpoonup z,\,hz\leftharpoonup z,\,ghz\leftharpoonup z$.
Then due to Eq.~\eqref{eq:r-action}, we actually only need to determine $z\rightharpoonup z,\,z\rightharpoonup gz,\,z\rightharpoonup hz,\,z\rightharpoonup ghz$. Finally, we utilize Eq.~\eqref{eq:MP1} to establish the relationship among $z\rightharpoonup z$, $z\rightharpoonup gz$, $z\rightharpoonup hz$ and $z\rightharpoonup ghz$, subsequently determine their specific values through the application of the method of undetermined coefficients.

\subsection{Classification of matched pairs of actions on $H_8$}\label{sec:cl}
\
\newline

\vspace{-.7em}
Suppose that $(H_8,\rightharpoonup,\leftharpoonup)$ is a matched pair of actions on $H_8$. Namely, Eqs.~\eqref{eq:MP1}--\eqref{eq:MP*} hold.
Then
\begin{align*}
&g\leftharpoonup g = g\rightharpoonup g,\\
&g\leftharpoonup h = (g\rightharpoonup h)gh,\\
&h\leftharpoonup g = (h\rightharpoonup g)gh,\\
&h\leftharpoonup h = h\rightharpoonup h.
\end{align*}
Since they are all non-trivial group-like elements, we have two situations to check:
\begin{enumerate}[(i)]
\item
$g\rightharpoonup h = h$,\quad $g\leftharpoonup h = g$;
\item
$g\rightharpoonup h = g$,\quad $g\leftharpoonup h = h$.
\end{enumerate}

First if $g\rightharpoonup h = h$ and $g\leftharpoonup h = g$, we must have $h\rightharpoonup g =g$ and $h\leftharpoonup g = h$. Otherwise,
$h\rightharpoonup g =h$ and $h\leftharpoonup g = g$, and then
$$h\leftharpoonup h = h\rightharpoonup h = h\rightharpoonup(h\rightharpoonup g)= g,$$
which implies that $g\leftharpoonup h = (h\leftharpoonup h)\leftharpoonup h=h$ as a contradiction. Also, $g\rightharpoonup g =h\rightharpoonup h =gh$ fails to hold, since
$$h\rightharpoonup h=h\rightharpoonup (g\rightharpoonup h)=gh\rightharpoonup h
=g\rightharpoonup (h\rightharpoonup h)=g\rightharpoonup (g\rightharpoonup g)=g$$
is a contradiction. As a result, we have the following 3 subcases:
$$(a)\ \left\{\begin{array}{l}
g\rightharpoonup g = g\leftharpoonup g = g,\\
g\rightharpoonup h = h\leftharpoonup g = h,\\
h\rightharpoonup g = g\leftharpoonup h = g,\\
h\rightharpoonup h = h\leftharpoonup h = h.
\end{array}\right.\quad (b)\ \left\{\begin{array}{l}
g\rightharpoonup g = g\leftharpoonup g = g,\\
g\rightharpoonup h = h\leftharpoonup g = h,\\
h\rightharpoonup g = g\leftharpoonup h = g,\\
h\rightharpoonup h = h\leftharpoonup h = gh.
\end{array}\right.\quad (c)\ \left\{\begin{array}{l}
g\rightharpoonup g = g\leftharpoonup g = gh,\\
g\rightharpoonup h = h\leftharpoonup g = h,\\
h\rightharpoonup g = g\leftharpoonup h = g,\\
h\rightharpoonup h = h\leftharpoonup h = h.
\end{array}\right.$$

On the other hand, if $g\rightharpoonup h = g$ and $g\leftharpoonup h = h$, then
$g\rightharpoonup g = g\rightharpoonup (g\rightharpoonup h)=h$ and
$h\leftharpoonup h = (g\leftharpoonup h)\leftharpoonup h = g$, so we have
$$(d)\ \left\{\begin{array}{l}
g\rightharpoonup g = g\leftharpoonup g = h,\\
g\rightharpoonup h = h\leftharpoonup g = g,\\
h\rightharpoonup g = g\leftharpoonup h = h,\\
h\rightharpoonup h = h\leftharpoonup h = g.
\end{array}\right.$$

By calculation, we find that no matched pair of actions on $H_8$ appears in cases (b)--(d). Take case (b) as an example. Let $x=h$ and $y=z$.
On the one hand,
\begin{align*}
&(x_1\rightharpoonup y_1) \otimes (x_2\leftharpoonup y_2)\\
&=\frac{1}{2}(h\rightharpoonup z\otimes h\leftharpoonup z+h\rightharpoonup gz\otimes h\leftharpoonup z+h\rightharpoonup z\otimes h\leftharpoonup hz-h\rightharpoonup gz\otimes h\leftharpoonup hz)\\
&=\frac{1}{2}(h\rightharpoonup z\otimes h\leftharpoonup z+g(h\rightharpoonup z)\otimes h\leftharpoonup z+h\rightharpoonup z\otimes gh\leftharpoonup z-g(h\rightharpoonup z)\otimes gh\leftharpoonup z).
\end{align*}
On the other hand,
\begin{align*}
&(x_2\rightharpoonup y_2) \otimes (x_1\leftharpoonup y_1)\\
&=\frac{1}{2}(h\rightharpoonup z\otimes h\leftharpoonup z+h\rightharpoonup z\otimes h\leftharpoonup gz+h\rightharpoonup hz\otimes h\leftharpoonup z-h\rightharpoonup hz\otimes h\leftharpoonup gz)\\
&=\frac{1}{2}(h\rightharpoonup z\otimes h\leftharpoonup z+h\rightharpoonup z\otimes h\leftharpoonup z).
\end{align*}
Then we have
$$(x_1\rightharpoonup y_1) \otimes (x_2\leftharpoonup y_2)-(x_2\rightharpoonup y_2) \otimes (x_1\leftharpoonup y_1)
=\frac{1}{2}(g-1)(h\rightharpoonup z)\otimes (1-g)h\leftharpoonup z\neq0,$$
so case (b) does not satisfy Eq.~\eqref{eq:MP5}.
Similarly, cases (c) and (d) do not satisfy Eq.~\eqref{eq:MP5} either.
Thus, we only need to consider case (a).

Suppose that a matched pair of actions $(H_8,\rightharpoonup, \leftharpoonup)$ satisfies case (a). Since
\begin{align*}
\Delta(z\rightharpoonup g)&=\frac{1}{2}(z \rightharpoonup g\otimes z \rightharpoonup g+
gz\rightharpoonup g\otimes z\rightharpoonup g+z\rightharpoonup g\otimes hz\rightharpoonup g-
gz\rightharpoonup g\otimes hz\rightharpoonup g)\\
&=\frac{1}{2}(z \rightharpoonup g\otimes z \rightharpoonup g+
zh\rightharpoonup g\otimes z\rightharpoonup g+z\rightharpoonup g\otimes zg\rightharpoonup g-
zh\rightharpoonup g\otimes zg\rightharpoonup g)\\
&=z\rightharpoonup g\otimes z\rightharpoonup g,
\end{align*}
we see that $z\rightharpoonup g \in G(H_8)$. Similarly, $z\rightharpoonup h \in G(H_8)$, $z\rightharpoonup gh \in G(H_8)$.

Also, we have $z\rightharpoonup g \neq 1$, for otherwise
\begin{equation*}
1\stackrel{\eqref{eq:MP2}}{=}z\rightharpoonup(z\rightharpoonup g)=z^{2}\rightharpoonup g=\frac{1}{2}(1+g+h-gh)\rightharpoonup g=g,
\end{equation*}which is a contradiction.
Similarly, $z\rightharpoonup h \neq 1$, $z\rightharpoonup gh \neq 1$.
As we know that
\begin{equation*}
z\rightharpoonup (z \rightharpoonup g)=g,\ z\rightharpoonup (z \rightharpoonup h)=h,\  z\rightharpoonup (z \rightharpoonup gh)=gh,
\end{equation*}
so we have 4 choices for the left action.
\begin{enumerate}[(i)]
\item $z\rightharpoonup g=g$, $z\rightharpoonup h=h$, $z\rightharpoonup gh=gh$.
\item $z\rightharpoonup g=g$, $z\rightharpoonup h=gh$, $z\rightharpoonup gh=h$.
\item $z\rightharpoonup g=gh$, $z\rightharpoonup h=h$, $z\rightharpoonup gh=g$.
\item $z\rightharpoonup g=h$, $z\rightharpoonup h=g$, $z\rightharpoonup gh=gh$.
\end{enumerate}

\medskip
Using the similar method, we can also obtain the following 4 choices for the right action.
\begin{enumerate}[(i)]
\item $g\leftharpoonup z=g$, $h\leftharpoonup z=h$, $gh\leftharpoonup z=gh$.
\item $g\leftharpoonup z=g$, $h\leftharpoonup z=gh$, $gh\leftharpoonup z=h$.
\item $g\leftharpoonup z=gh$, $h\leftharpoonup z=h$, $gh\leftharpoonup z=g$.
\item $g\leftharpoonup z=h$, $h\leftharpoonup z=g$, $gh\leftharpoonup z=gh$.
\end{enumerate}
Therefore, we conclude that in order to find all the matched pairs of actions on $H_8$, there are 16 situations to be considered. We will only discuss in detail two of them which really contribute to the classification.

\subsubsection{Situation 1}
The first situation we discuss is as follows.
$$z\rightharpoonup g=g,\ z\rightharpoonup h=h,\ z\rightharpoonup gh=gh,\ g\leftharpoonup z=g,\ h\leftharpoonup z=h,\ gh\leftharpoonup z=gh.$$
Then we have
\begin{align*}
g\rightharpoonup z&\stackrel{\eqref{eq:l-action}}=\frac{1}{2}(gzS(g\leftharpoonup z)+ggzS(g\leftharpoonup z)+gzS(g\leftharpoonup hz)-ggzS(g\leftharpoonup hz))\\
&=\frac{1}{2}(gz(g\leftharpoonup z)+z(g\leftharpoonup z)+gz(g\leftharpoonup z)-z(g\leftharpoonup z))\\
&=gz(g\leftharpoonup z)\\
&=ghz.
\end{align*}
Similarly, $h\rightharpoonup z=hz(h\leftharpoonup z)=ghz$, $gh\rightharpoonup z=ghz(gh\leftharpoonup z)=z$. Also, we calculate that
\begin{align*}
z\rightharpoonup gz&\stackrel{\eqref{eq:MP1}}=\frac{1}{2}((z\rightharpoonup g)((z\leftharpoonup g)\rightharpoonup z)+(gz\rightharpoonup g)((z\leftharpoonup g)\rightharpoonup z)\\
&\quad +(z\rightharpoonup g)((hz\leftharpoonup g)\rightharpoonup z)-(gz\rightharpoonup g)((hz\leftharpoonup g)\rightharpoonup z))\\
&=(z\rightharpoonup g)((z\leftharpoonup g)\rightharpoonup z)\\
&\stackrel{\eqref{eq:r-action}}=g(ghz\rightharpoonup z)\\
&=g(z\rightharpoonup z).
\end{align*}
Similarly, $z\rightharpoonup hz=h(z\rightharpoonup z)$, $z\rightharpoonup ghz=gh(z\rightharpoonup z)$. Correspondingly, we see that
$$gz\rightharpoonup z =zh\rightharpoonup z =z\rightharpoonup (h\rightharpoonup z) = z\rightharpoonup ghz =gh(z\rightharpoonup z),$$
$hz\rightharpoonup z=gh(z\rightharpoonup z)$ and $ghz\rightharpoonup z=z\rightharpoonup z$. Then the other nine elements
$$Gz\rightharpoonup G'z,\quad G,G'\in\{g,h,gh\}$$
can be checked by Eq.~\eqref{eq:MP1}.
In summary, we have
\begin{equation}\label{eq:z-z}
\begin{cases}
gz\rightharpoonup ghz=hz\rightharpoonup ghz=ghz\rightharpoonup z=z\rightharpoonup z,\\
z\rightharpoonup gz=gz\rightharpoonup hz=hz\rightharpoonup hz=ghz\rightharpoonup gz=g(z\rightharpoonup z),\\
z\rightharpoonup hz=gz\rightharpoonup gz=hz\rightharpoonup gz=ghz\rightharpoonup hz=h(z\rightharpoonup z),\\
z\rightharpoonup ghz=
gz\rightharpoonup z=hz\rightharpoonup z=ghz\rightharpoonup ghz=gh(z\rightharpoonup z).
\end{cases}
\end{equation}
Hence, the left action $\rightharpoonup$ will be figured out, whenever $z\rightharpoonup z$ is determined.

To compute $z\rightharpoonup z$, we set
$$z\rightharpoonup z=a_1+a_2g+a_3h+a_4gh+a_5z+a_6gz+a_7hz+a_8ghz.$$
On the one hand,
\begin{equation*}
z^2\rightharpoonup z=\frac{1}{2}(1+g+h-gh)\rightharpoonup z=\frac{1}{2}(z+ghz+ghz-z)=ghz.
\end{equation*}
On the other hand,
\begin{align*}
&z^2\rightharpoonup z=z\rightharpoonup (z\rightharpoonup z)\\
&=z\rightharpoonup (a_1+a_2g+a_3h+a_4gh+a_5z+a_6gz+a_7hz+a_8ghz)\\
&=a_1+a_2g+a_3h+a_4gh+a_5(a_1+a_2g+a_3h+a_4gh+a_5z+a_6gz+a_7hz+a_8ghz)\\
&\quad +a_6(a_1g+a_2+a_3gh+a_4h+a_5gz+a_6z+a_7ghz+a_8hz)\\
&\quad +a_7(a_1h+a_2gh+a_3+a_4g+a_5hz+a_6ghz+a_7z+a_8gz)\\
&\quad +a_8(a_1gh+a_2h+a_3g+a_4+a_5ghz+a_6hz+a_7gz+a_8z).
\end{align*}
By the comparison of coefficients, we have
\begin{equation}\label{fangchengzu1}
\begin{cases}
a_1+a_5a_1+a_6a_2+a_7a_3+a_8a_4=0,\\
a_2+a_5a_2+a_6a_1+a_7a_4+a_8a_3=0,\\
a_3+a_5a_3+a_6a_4+a_7a_1+a_8a_2=0,\\
a_4+a_5a_4+a_6a_3+a_7a_2+a_8a_1=0,\\
a_5^2+a_6^2+a_7^2+a_8^2=0,\\
a_5a_6+a_6a_5+a_7a_8+a_8a_7=0,\\
a_5a_7+a_6a_8+a_7a_5+a_8a_6=0,\\
a_5a_8+a_6a_7+a_7a_6+a_8a_5=1,
\end{cases}
\end{equation}
which clearly implies that
$$(a_1+a_2+a_3+a_4)(a_5+a_6+a_7+a_8+1)=0,\ (a_5+a_6+a_7+a_8)^2=1.$$
Also, the left module coalgebra action $\rightharpoonup$ on $H_8$ satisfies $\varepsilon(z\rightharpoonup z)=\varepsilon(z)\varepsilon(z)$. That is, $a_1+a_2+a_3+a_4+a_5+a_6+a_7+a_8=1$.
Hence, we  get two choices.
\begin{enumerate}[(1)]
\item
$a_1+a_2+a_3+a_4=2,\ a_5+a_6+a_7+a_8=-1$,
\item
$a_1+a_2+a_3+a_4=0,\ a_5+a_6+a_7+a_8=1$.
\end{enumerate}

But none of them is enough to determine $z\rightharpoonup z$, so we consider
the coproduct $\Delta(z\rightharpoonup z)$ to get more constraints.

On the one hand, $\Delta(z\rightharpoonup z)\in \Delta(a_5z+a_6gz+a_7hz+a_8ghz)+\bk G(H_8)\otimes \bk G(H_8)$ and
\begin{align*}
\Delta&(a_5z+a_6gz+a_7hz+a_8ghz)
=a_5\Delta(z)+a_6(g\otimes g)\Delta(z)+a_7(h\otimes h)\Delta(z)+a_8(gh\otimes gh)\Delta(z)\\
&=\frac{1}{2}a_5(z\otimes z+gz\otimes z+z\otimes hz-gz\otimes hz)+\frac{1}{2}a_6(gz\otimes gz+z\otimes gz+gz\otimes ghz-z\otimes ghz)\\
&\quad +\frac{1}{2}a_7(hz\otimes hz+ghz\otimes hz+hz\otimes z-ghz\otimes z)+\frac{1}{2}a_8(ghz\otimes ghz+hz\otimes ghz+ghz\otimes gz-hz\otimes gz).
\end{align*}
On the other hand,
\begin{equation}\label{delta}
\begin{aligned}
&\Delta(z\rightharpoonup z)=\frac{1}{4}(z\rightharpoonup z\otimes z\rightharpoonup z+z\rightharpoonup gz\otimes z\rightharpoonup z+z\rightharpoonup z\otimes z\rightharpoonup hz-z\rightharpoonup gz\otimes z\rightharpoonup hz\\
&\quad+gz\rightharpoonup z\otimes z\rightharpoonup z+gz\rightharpoonup gz\otimes z\rightharpoonup z+gz\rightharpoonup z\otimes z\rightharpoonup hz-gz\rightharpoonup gz\otimes z\rightharpoonup hz\\
&\quad+z\rightharpoonup z\otimes hz\rightharpoonup z+z\rightharpoonup gz\otimes hz\rightharpoonup z+z\rightharpoonup z\otimes hz\rightharpoonup hz-z\rightharpoonup gz\otimes hz\rightharpoonup hz\\
&\quad-gz\rightharpoonup z\otimes hz\rightharpoonup z-gz\rightharpoonup gz\otimes hz\rightharpoonup z-gz\rightharpoonup z\otimes hz\rightharpoonup hz+gz\rightharpoonup gz\otimes hz\rightharpoonup hz)\\
&\stackrel{\eqref{eq:z-z}}=
\frac{1}{4}(1\otimes1+g\otimes1+1\otimes h -g\otimes h+gh\otimes1 +h\otimes1 + gh\otimes h -h\otimes h+1\otimes gh \\
&\quad+g\otimes gh+1\otimes g -g\otimes g -gh\otimes gh -h\otimes gh -gh\otimes g +h\otimes g)(z\rightharpoonup z\otimes z\rightharpoonup z).
\end{aligned}
\end{equation}

By comparing coefficients of $Gz\otimes G'z$, where $G,G'\in\{1, g,h,gh\}$, in the above two coproducts, we obtain the following 16 equations,
{\footnotesize\begin{equation}\label{fangchengzu2}
\begin{cases}
a_5^2-a_6^2-a_7^2-a_8^2+2a_5a_6+2a_5a_7+2a_5a_8=2a_5,\\
a_5^2+a_6^2-a_7^2+a_8^2+2a_5a_7+2a_6a_7-2a_7a_8=2a_6,\\
a_5^2+a_6^2+a_7^2-a_8^2+2a_5a_8-2a_6a_8+2a_7a_8=2a_5,\\
a_5^2-a_6^2+a_7^2+a_8^2+2a_5a_6-2a_6a_7+2a_6a_8=-2a_6,\\
a_5^2+a_6^2+a_7^2-a_8^2+2a_5a_8+2a_6a_8-2a_7a_8=2a_5,\\
-a_5^2+a_6^2-a_7^2-a_8^2+2a_5a_6+2a_6a_7+2a_6a_8=2a_6,\\
-a_5^2+a_6^2+a_7^2+a_8^2+2a_5a_6+2a_5a_7-2a_5a_8=-2a_5,\\
a_5^2+a_6^2-a_7^2+a_8^2-2a_5a_7+2a_6a_7+2a_7a_8=2a_6,\\
a_5^2-a_6^2+a_7^2+a_8^2+2a_5a_6+2a_6a_7-2a_6a_8=2a_7,\\
a_5^2+a_6^2+a_7^2-a_8^2-2a_5a_8+2a_6a_8+2a_7a_8=-2a_8,\\
-a_5^2-a_6^2+a_7^2-a_8^2+2a_5a_7+2a_6a_7+2a_7a_8=2a_7,\\
-a_5^2+a_6^2+a_7^2+a_8^2-2a_5a_6+2a_5a_7+2a_5a_8=2a_8,\\
a_5^2+a_6^2-a_7^2+a_8^2+2a_5a_7-2a_6a_7+2a_7a_8=-2a_7,\\
-a_5^2+a_6^2+a_7^2+a_8^2+2a_5a_6-2a_5a_7+2a_5a_8=2a_8,\\
a_5^2-a_6^2+a_7^2+a_8^2-2a_5a_6+2a_6a_7+2a_6a_8=2a_7,\\
-a_5^2-a_6^2-a_7^2+a_8^2+2a_5a_8+2a_6a_8+2a_7a_8=2a_8.
\end{cases}
\end{equation}}
By dividing them into 8 pairs with terms cancelled, we get that
$a_5a_6=a_5a_7=a_6a_8=a_7a_8=0$.
Together with Eq.~\eqref{fangchengzu1}, it means that $a_5=a_8=0$, $a_6a_7=\frac{1}{2}$ or $a_6=a_7=0$, $a_5a_8=\frac{1}{2}$.

Now we come back to discuss the two choices obtained from Eq.~\eqref{fangchengzu1}. When $a_1+a_2+a_3+a_4=2,\,a_5+a_6+a_7+a_8=-1$, one more constraint $a_5=a_8=0$, $a_6a_7=\frac{1}{2}$ is incompatible with the second equation in \eqref{fangchengzu2},
while one more constraint $a_6=a_7=0$, $a_5a_8=\frac{1}{2}$ contradicts to the first one in \eqref{fangchengzu2}. When $a_1+a_2+a_3+a_4=0,\,a_5+a_6+a_7+a_8=1$, we can obtain the following four solutions.

\smallskip
\textbf{Case 1.1}
\quad$a_5=0,\,a_6=\frac{1-\sqrt{-1}}{2},\,a_7=\frac{1+\sqrt{-1}}{2},\,a_8=0.$

\smallskip
\textbf{Case 1.2}
\quad$a_5=0,\,a_6=\frac{1+\sqrt{-1}}{2},\,a_7=\frac{1-\sqrt{-1}}{2},\,a_8=0.$

\smallskip
\textbf{Case 1.3}
\quad$a_5=\frac{1-\sqrt{-1}}{2},\,a_6=0,\,a_7=0,\,a_8=\frac{1+\sqrt{-1}}{2}.$

\smallskip
\textbf{Case 1.4}
\quad$a_5=\frac{1+\sqrt{-1}}{2},\,a_6=0,\,a_7=0,\,a_8=\frac{1-\sqrt{-1}}{2}.$

\smallskip
After determining $a_5$, $a_6$, $a_7$, $a_8$, we will proceed to calculate $a_1$, $a_2$, $a_3$, $a_4$. Taking Case 1.1 as an example, we have
$$z\rightharpoonup z=a_1+a_2g+a_3h+a_4gh+\frac{1-\sqrt{-1}}{2}gz+\frac{1+\sqrt{-1}}{2}hz,$$
and then the coproduct
$$\Delta(z\rightharpoonup z)=a_11\otimes 1+a_2g\otimes g+a_3h\otimes h+a_4gh\otimes gh+\frac{1-\sqrt{-1}}{2}(g\otimes g)\Delta(z)+\frac{1+\sqrt{-1}}{2}(h\otimes h)\Delta(z)$$
does not contain any term of which the first tensor factor is a group-like element and the second one is $z$, $gz$, $hz$ or $ghz$. On the other hand, Eq.~\eqref{delta} tells us that the sum of all these terms in $\Delta(z\rightharpoonup z)$ is
\begin{align*}
&\frac{1}{4}((1+g+h+gh)\otimes1 +(1-g+h-gh)\otimes g +(1-g-h+gh)\otimes h  +(1+g-h-gh)\otimes gh)\\
&\quad ((a_1+a_2g+a_3h+a_4gh)\otimes (\frac{1-\sqrt{-1}}{2}gz+\frac{1+\sqrt{-1}}{2}hz)),
\end{align*}
and it vanishes if and only if
$$a_1=a_2=a_3=a_4=0.$$
In the remaining 3 cases, we can similarly obtain that $a_1=a_2=a_3=a_4=0$ as in Case 1.1.
So far all cases for $z\rightharpoonup z$ have been considered, and we list below the multiplication tables for the left actions $\rightharpoonup$, which have been verified as left module coalgebra actions.
\begin{table}[H]
\centering
\caption{The left action on $H_8$ in Case 1.1}\label{H8}
\renewcommand{\arraystretch}{1.4}
\begin{tabular}{|c|c|c|c|c|c|c|c|c|}
\hline
$\rightharpoonup$&1&$g$&$h$&$gh$&$z$&$gz$&$hz$&$ghz$\\
\hline
1&1&$g$&$h$&$gh$&$z$&$gz$&$hz$&$ghz$\\
\hline
$g$&1&$g$&$h$&$gh$&$ghz$&$hz$&$gz$&$z$\\
\hline
$h$&$1$&$g$&$h$&$gh$&$ghz$&$hz$&$gz$&$z$\\
\hline
$gh$&1&$g$&$h$&$gh$&$z$&$gz$&$hz$&$ghz$\\
\hline
$z$&1&$g$&$h$&$gh$&$\frac{1-\sqrt{-1}}{2}gz$+$\frac{1+\sqrt{-1}}{2}hz$&$\frac{1-\sqrt{-1}}{2}z$+$\frac{1+\sqrt{-1}}{2}ghz$&$\frac{1-\sqrt{-1}}{2}ghz$+$\frac{1+\sqrt{-1}}{2}z$&$\frac{1-\sqrt{-1}}{2}hz$+$\frac{1+\sqrt{-1}}{2}gz$\\
\hline
$gz$&1&$g$&$h$&$gh$&$\frac{1-\sqrt{-1}}{2}hz$+$\frac{1+\sqrt{-1}}{2}gz$&$\frac{1-\sqrt{-1}}{2}ghz$+$\frac{1+\sqrt{-1}}{2}z$&$\frac{1-\sqrt{-1}}{2}z$+$\frac{1+\sqrt{-1}}{2}ghz$&$\frac{1-\sqrt{-1}}{2}gz$+$\frac{1+\sqrt{-1}}{2}hz$\\
\hline
$hz$&$1$&$g$&$h$&$gh$&$\frac{1-\sqrt{-1}}{2}hz$+$\frac{1+\sqrt{-1}}{2}gz$&$\frac{1-\sqrt{-1}}{2}ghz$+$\frac{1+\sqrt{-1}}{2}z$&$\frac{1-\sqrt{-1}}{2}z$+$\frac{1+\sqrt{-1}}{2}ghz$&$\frac{1-\sqrt{-1}}{2}gz$+$\frac{1+\sqrt{-1}}{2}hz$\\
\hline
$ghz$&1&$g$&$h$&$gh$&$\frac{1-\sqrt{-1}}{2}gz$+$\frac{1+\sqrt{-1}}{2}hz$&$\frac{1-\sqrt{-1}}{2}z$+$\frac{1+\sqrt{-1}}{2}ghz$&$\frac{1-\sqrt{-1}}{2}ghz$+$\frac{1+\sqrt{-1}}{2}z$&$\frac{1-\sqrt{-1}}{2}hz$+$\frac{1+\sqrt{-1}}{2}gz$\\
\hline
\end{tabular}

\bigskip
\centering
\caption{The left action on $H_8$ in Case 1.2}\label{H8'}
\renewcommand{\arraystretch}{1.4}
\begin{tabular}{|c|c|c|c|c|c|c|c|c|}
\hline
$\rightharpoonup$&1&$g$&$h$&$gh$&$z$&$gz$&$hz$&$ghz$\\
\hline
1&1&$g$&$h$&$gh$&$z$&$gz$&$hz$&$ghz$\\
\hline
$g$&1&$g$&$h$&$gh$&$ghz$&$hz$&$gz$&$z$\\
\hline
$h$&$1$&$g$&$h$&$gh$&$ghz$&$hz$&$gz$&$z$\\
\hline
$gh$&1&$g$&$h$&$gh$&$z$&$gz$&$hz$&$ghz$\\
\hline
$z$&1&$g$&$h$&$gh$&$\frac{1+\sqrt{-1}}{2}gz$+$\frac{1-\sqrt{-1}}{2}hz$&$\frac{1+\sqrt{-1}}{2}z$+$\frac{1-\sqrt{-1}}{2}ghz$&$\frac{1+\sqrt{-1}}{2}ghz$+$\frac{1-\sqrt{-1}}{2}z$&$\frac{1+\sqrt{-1}}{2}hz$+$\frac{1-\sqrt{-1}}{2}gz$\\
\hline
$gz$&1&$g$&$h$&$gh$&$\frac{1+\sqrt{-1}}{2}hz$+$\frac{1-\sqrt{-1}}{2}gz$&$\frac{1+\sqrt{-1}}{2}ghz$+$\frac{1-\sqrt{-1}}{2}z$&$\frac{1+\sqrt{-1}}{2}z$+$\frac{1-\sqrt{-1}}{2}ghz$&$\frac{1+\sqrt{-1}}{2}gz$+$\frac{1-\sqrt{-1}}{2}hz$\\
\hline
$hz$&$1$&$g$&$h$&$gh$&$\frac{1+\sqrt{-1}}{2}hz$+$\frac{1-\sqrt{-1}}{2}gz$&$\frac{1+\sqrt{-1}}{2}ghz$+$\frac{1-\sqrt{-1}}{2}z$&$\frac{1+\sqrt{-1}}{2}z$+$\frac{1-\sqrt{-1}}{2}ghz$&$\frac{1+\sqrt{-1}}{2}gz$+$\frac{1-\sqrt{-1}}{2}hz$\\
\hline
$ghz$&1&$g$&$h$&$gh$&$\frac{1+\sqrt{-1}}{2}gz$+$\frac{1-\sqrt{-1}}{2}hz$&$\frac{1+\sqrt{-1}}{2}z$+$\frac{1-\sqrt{-1}}{2}ghz$&$\frac{1+\sqrt{-1}}{2}ghz$+$\frac{1-\sqrt{-1}}{2}z$&$\frac{1+\sqrt{-1}}{2}hz$+$\frac{1-\sqrt{-1}}{2}gz$\\
\hline
\end{tabular}

\bigskip
\centering
\caption{The left action on $H_8$ in Case 1.3}\label{H8''}
\renewcommand{\arraystretch}{1.4}
\begin{tabular}{|c|c|c|c|c|c|c|c|c|}
\hline
$\rightharpoonup$&1&$g$&$h$&$gh$&$z$&$gz$&$hz$&$ghz$\\
\hline
1&1&$g$&$h$&$gh$&$z$&$gz$&$hz$&$ghz$\\
\hline
$g$&1&$g$&$h$&$gh$&$ghz$&$hz$&$gz$&$z$\\
\hline
$h$&$1$&$g$&$h$&$gh$&$ghz$&$hz$&$gz$&$z$\\
\hline
$gh$&1&$g$&$h$&$gh$&$z$&$gz$&$hz$&$ghz$\\
\hline
$z$&1&$g$&$h$&$gh$&$\frac{1-\sqrt{-1}}{2}z$+$\frac{1+\sqrt{-1}}{2}ghz$&$\frac{1-\sqrt{-1}}{2}gz$+$\frac{1+\sqrt{-1}}{2}hz$&$\frac{1-\sqrt{-1}}{2}hz$+$\frac{1+\sqrt{-1}}{2}gz$&$\frac{1-\sqrt{-1}}{2}ghz$+$\frac{1+\sqrt{-1}}{2}z$\\
\hline
$gz$&1&$g$&$h$&$gh$&$\frac{1-\sqrt{-1}}{2}ghz$+$\frac{1+\sqrt{-1}}{2}z$&$\frac{1-\sqrt{-1}}{2}hz$+$\frac{1+\sqrt{-1}}{2}gz$&$\frac{1-\sqrt{-1}}{2}gz$+$\frac{1+\sqrt{-1}}{2}hz$&$\frac{1-\sqrt{-1}}{2}z$+$\frac{1+\sqrt{-1}}{2}ghz$\\
\hline
$hz$&$1$&$g$&$h$&$gh$&$\frac{1-\sqrt{-1}}{2}ghz$+$\frac{1+\sqrt{-1}}{2}z$&$\frac{1-\sqrt{-1}}{2}hz$+$\frac{1+\sqrt{-1}}{2}gz$&$\frac{1-\sqrt{-1}}{2}gz$+$\frac{1+\sqrt{-1}}{2}hz$&$\frac{1-\sqrt{-1}}{2}z$+$\frac{1+\sqrt{-1}}{2}ghz$\\
\hline
$ghz$&1&$g$&$h$&$gh$&$\frac{1-\sqrt{-1}}{2}z$+$\frac{1+\sqrt{-1}}{2}ghz$&$\frac{1-\sqrt{-1}}{2}gz$+$\frac{1+\sqrt{-1}}{2}hz$&$\frac{1-\sqrt{-1}}{2}hz$+$\frac{1+\sqrt{-1}}{2}gz$&$\frac{1-\sqrt{-1}}{2}ghz$+$\frac{1+\sqrt{-1}}{2}z$\\
\hline
\end{tabular}
\end{table}
\begin{table}[H]
\centering
\caption{The left action on $H_8$ in Case 1.4}\label{H8'''}
\renewcommand{\arraystretch}{1.4}
\begin{tabular}{|c|c|c|c|c|c|c|c|c|}
\hline
$\rightharpoonup$&1&$g$&$h$&$gh$&$z$&$gz$&$hz$&$ghz$\\
\hline
1&1&$g$&$h$&$gh$&$z$&$gz$&$hz$&$ghz$\\
\hline
$g$&1&$g$&$h$&$gh$&$ghz$&$hz$&$gz$&$z$\\
\hline
$h$&$1$&$g$&$h$&$gh$&$ghz$&$hz$&$gz$&$z$\\
\hline
$gh$&1&$g$&$h$&$gh$&$z$&$gz$&$hz$&$ghz$\\
\hline
$z$&1&$g$&$h$&$gh$&$\frac{1+\sqrt{-1}}{2}z$+$\frac{1-\sqrt{-1}}{2}ghz$&$\frac{1+\sqrt{-1}}{2}gz$+$\frac{1-\sqrt{-1}}{2}hz$&$\frac{1+\sqrt{-1}}{2}hz$+$\frac{1-\sqrt{-1}}{2}gz$&$\frac{1+\sqrt{-1}}{2}ghz$+$\frac{1-\sqrt{-1}}{2}z$\\
\hline
$gz$&1&$g$&$h$&$gh$&$\frac{1+\sqrt{-1}}{2}ghz$+$\frac{1-\sqrt{-1}}{2}z$&$\frac{1+\sqrt{-1}}{2}hz$+$\frac{1-\sqrt{-1}}{2}gz$&$\frac{1+\sqrt{-1}}{2}gz$+$\frac{1-\sqrt{-1}}{2}hz$&$\frac{1+\sqrt{-1}}{2}z$+$\frac{1-\sqrt{-1}}{2}ghz$\\
\hline
$hz$&$1$&$g$&$h$&$gh$&$\frac{1+\sqrt{-1}}{2}ghz$+$\frac{1-\sqrt{-1}}{2}z$&$\frac{1+\sqrt{-1}}{2}hz$+$\frac{1-\sqrt{-1}}{2}gz$&$\frac{1+\sqrt{-1}}{2}gz$+$\frac{1-\sqrt{-1}}{2}hz$&$\frac{1+\sqrt{-1}}{2}z$+$\frac{1-\sqrt{-1}}{2}ghz$\\
\hline
$ghz$&1&$g$&$h$&$gh$&$\frac{1+\sqrt{-1}}{2}z$+$\frac{1-\sqrt{-1}}{2}ghz$&$\frac{1+\sqrt{-1}}{2}gz$+$\frac{1-\sqrt{-1}}{2}hz$&$\frac{1+\sqrt{-1}}{2}hz$+$\frac{1-\sqrt{-1}}{2}gz$&$\frac{1+\sqrt{-1}}{2}ghz$+$\frac{1-\sqrt{-1}}{2}z$\\
\hline
\end{tabular}
\end{table}

\subsubsection{Situation 2}
The second situation we discuss is as follows.
$$z\rightharpoonup g=h,\ z\rightharpoonup h=g,\ z\rightharpoonup gh=gh,\ g\leftharpoonup z=h,\ h\leftharpoonup z=g,\ gh\leftharpoonup z=gh.$$
Similar to the Situation 1, using Eqs.~\eqref{eq:l-action}, \eqref{eq:MP1} and \eqref{eq:r-action}, we compute that
\begin{equation}\label{eq:z-z'}
\begin{cases}
g\rightharpoonup z=h\rightharpoonup z=gh\rightharpoonup z=z,\\
gz\rightharpoonup z=hz\rightharpoonup z=ghz\rightharpoonup z=z\rightharpoonup z,\\
z\rightharpoonup hz=gz\rightharpoonup hz=hz\rightharpoonup hz=ghz\rightharpoonup hz=g(z\rightharpoonup z),\\
z\rightharpoonup gz=gz\rightharpoonup gz=hz\rightharpoonup gz=ghz\rightharpoonup gz=h(z\rightharpoonup z),\\
z\rightharpoonup ghz=gz\rightharpoonup ghz=hz\rightharpoonup ghz=ghz\rightharpoonup ghz=gh(z\rightharpoonup z).
\end{cases}
\end{equation}

So again it is crucial to determine $z\rightharpoonup z$. Set
$$z\rightharpoonup z=a_1+a_2g+a_3h+a_4gh+a_5z+a_6gz+a_7hz+a_8ghz.$$
On the one hand, $$z^2\rightharpoonup z=\frac{1}{2}(1+g+h-gh)\rightharpoonup z=\frac{1}{2}(z+z+z-z)=z.$$
On the other hand,
\begin{equation}\label{xishu}
\begin{aligned}
&z^2\rightharpoonup z=z\rightharpoonup (z\rightharpoonup z)\\
&=a_1+a_2h+a_3g+a_4gh+a_5(a_1+a_2g+a_3h+a_4gh+a_5z+a_6gz+a_7hz+a_8ghz)\\
&\quad+a_6(a_1h+a_2gh+a_3+a_4g+a_5hz+a_6ghz+a_7z+a_8gz)\\
&\quad+a_7(a_1g+a_2+a_3gh+a_4h+a_5gz+a_6z+a_7ghz+a_8hz)\\
&\quad+a_8(a_1gh+a_2h+a_3g+a_4+a_5ghz+a_6hz+a_7gz+a_8z).
\end{aligned}
\end{equation}
By the comparison of coefficients, we have
\begin{equation}\label{fangchengzu3}
\begin{cases}
a_1+a_5a_1+a_6a_3+a_7a_2+a_8a_4=0,\\
a_3+a_5a_2+a_6a_4+a_7a_1+a_8a_3=0,\\
a_2+a_5a_3+a_6a_1+a_7a_4+a_8a_2=0,\\
a_4+a_5a_4+a_6a_2+a_7a_3+a_8a_1=0,\\
a_5^2+a_6a_7+a_7a_6+a_8^2=1,\\
a_5a_6+a_6a_8+a_7a_5+a_8a_7=0,\\
a_5a_7+a_6a_5+a_7a_8+a_8a_6=0,\\
a_5a_8+a_6^2+a_7^2+a_8a_5=0.\\
\end{cases}
\end{equation}
As $a_1+a_2+a_3+a_4+a_5+a_6+a_7+a_8=1$, we also have
two choices in this situation.
\begin{enumerate}[(1)]
\item
$a_1+a_2+a_3+a_4=2,\ a_5+a_6+a_7+a_8=-1$,
\item
$a_1+a_2+a_3+a_4=0,\ a_5+a_6+a_7+a_8=1$.
\end{enumerate}
Analogous to Situation 1, we will first calculate $a_5$, $a_6$, $a_7$, $a_8$.
On the one hand,
\begin{align*}
&\Delta(a_5z+a_6gz+a_7hz+a_8ghz)\\
&=\frac{1}{2}a_5(z\otimes z+gz\otimes z+z\otimes hz-gz\otimes hz)+\frac{1}{2}a_6(gz\otimes gz+z\otimes gz+gz\otimes ghz-z\otimes ghz)\\
&\quad+\frac{1}{2}a_7(hz\otimes hz+ghz\otimes hz+hz\otimes z-ghz\otimes z)+\frac{1}{2}a_8(ghz\otimes ghz+hz\otimes ghz+ghz\otimes gz-hz\otimes gz).\\
\end{align*}
On the other hand,
$$\Delta(z\rightharpoonup z)\stackrel{\eqref{eq:z-z'}}=\frac{1}{2}(1\otimes1+h\otimes1+1\otimes g-h\otimes g)(z\rightharpoonup z\otimes z\rightharpoonup z).$$
By comparing coefficients of $Gz\otimes G'z$, where $G,G'\in\{1, g,h,gh\}$, in the above two coproducts, we get $a_5^2=a_6^2=a_7^2=a_8^2$. Then using
Eq.~\eqref{fangchengzu3}, we can obtain that $a_5=a_6=a_7=-a_8$ or $-a_5=a_6=a_7=a_8$.

Now for Choice (1) $a_1+a_2+a_3+a_4=2,\,a_5+a_6+a_7+a_8=-1$, no solution is found.

For Choice (2) $a_1+a_2+a_3+a_4=0,\,a_5+a_6+a_7+a_8=1$, we obtain the following two solutions.

\smallskip
\textbf{Case 2.1}
\quad$a_5=\frac{1}{2},\,a_6=\frac{1}{2},\,a_7=\frac{1}{2},\,a_8=-\frac{1}{2}.$

\smallskip
\textbf{Case 2.2}
\quad$a_5=-\frac{1}{2},\,a_6=\frac{1}{2},\,a_7=\frac{1}{2},\,a_8=\frac{1}{2}.$

\smallskip
In these two cases, we can calculate that $a_1=a_2=a_3=a_4=0$
using the same method as Situation 1 shown.  The final results are shown as follows, where $z^{3}=\frac{1}{2}(z+gz+hz-ghz)$.

\begin{table}[H]
\centering
\begin{minipage}{0.45\textwidth}
\centering
\caption{The left action on  \hspace*{4em} $H_8$ in Case 2.1}\label{H8-}
\renewcommand{\arraystretch}{1.2}
\begin{tabular}{|c|c|c|c|c|c|c|c|c|}
\hline
$\rightharpoonup$&1&$g$&$h$&$gh$&$z$&$gz$&$hz$&$ghz$\\
\hline
1&1&$g$&$h$&$gh$&$z$&$gz$&$hz$&$ghz$\\
\hline
$g$&1&$g$&$h$&$gh$&$z$&$gz$&$hz$&$ghz$\\
\hline
$h$&1&$g$&$h$&$gh$&$z$&$gz$&$hz$&$ghz$\\
\hline
$gh$&1&$g$&$h$&$gh$&$z$&$gz$&$hz$&$ghz$\\
\hline
$z$&1&$h$&$g$&$gh$&$z^{3}$&$hz^{3}$&$gz^{3}$&$ghz^{3}$\\
\hline
$gz$&1&$h$&$g$&$gh$&$z^{3}$&$hz^{3}$&$gz^{3}$&$ghz^{3}$\\
\hline
$hz$&1&$h$&$g$&$gh$&$z^{3}$&$hz^{3}$&$gz^{3}$&$ghz^{3}$\\
\hline
$ghz$&1&$h$&$g$&$gh$&$z^{3}$&$hz^{3}$&$gz^{3}$&$ghz^{3}$\\
\hline
\end{tabular}
\end{minipage}
\hspace{0.05\textwidth} 
\begin{minipage}{0.45\textwidth}
\centering
\caption{The left action on \hspace*{4em} $H_8$ in Case 2.2}\label{H8--}
\renewcommand{\arraystretch}{1.2}
\begin{tabular}{|c|c|c|c|c|c|c|c|c|}
\hline
$\rightharpoonup$&1&$g$&$h$&$gh$&$z$&$gz$&$hz$&$ghz$\\
\hline
1&1&$g$&$h$&$gh$&$z$&$gz$&$hz$&$ghz$\\
\hline
$g$&1&$g$&$h$&$gh$&$z$&$gz$&$hz$&$ghz$\\
\hline
$h$&1&$g$&$h$&$gh$&$z$&$gz$&$hz$&$ghz$\\
\hline
$gh$&1&$g$&$h$&$gh$&$z$&$gz$&$hz$&$ghz$\\
\hline
$z$&1&$h$&$g$&$gh$&$ghz^{3}$&$gz^{3}$&$hz^{3}$&$z^{3}$\\
\hline
$gz$&1&$h$&$g$&$gh$&$ghz^{3}$&$gz^{3}$&$hz^{3}$&$z^{3}$\\
\hline
$hz$&1&$h$&$g$&$gh$&$ghz^{3}$&$gz^{3}$&$hz^{3}$&$z^{3}$\\
\hline
$ghz$&1&$h$&$g$&$gh$&$ghz^{3}$&$gz^{3}$&$hz^{3}$&$z^{3}$\\
\hline
\end{tabular}
\end{minipage}
\end{table}

\begin{remark}
For the remaining 14 situations, we discuss them similarly to the above two cases. However, when we determine
$z\rightharpoonup z$ by the method of undetermined coefficients, no solution can be found in each of these situations.
So they do not give any matched pair of actions on $H_8$.

Take the following situation as an example: $$z\rightharpoonup g=h,\ z\rightharpoonup h=g,\ z\rightharpoonup gh=gh,\ g\leftharpoonup z=g,\ h\leftharpoonup z=gh,\ gh\leftharpoonup z=h.$$
Similar to Situation 2, using Eqs.~\eqref{eq:l-action} and \eqref{eq:MP1}, we obtain that
\begin{align*}
&g\rightharpoonup z=gz(g\leftharpoonup z)=ghz,\,h\rightharpoonup z=hz(h\leftharpoonup z)=gz,\,gh\rightharpoonup z=ghz(gh\leftharpoonup z)=hz,\\
&z\rightharpoonup gz=h(z\rightharpoonup z),\,z\rightharpoonup hz=g(z\rightharpoonup z),\,z\rightharpoonup ghz=gh(z\rightharpoonup z).
\end{align*}

Set $z\rightharpoonup z=a_1+a_2g+a_3h+a_4gh+a_5z+a_6gz+a_7hz+a_8ghz$.
Combining
$$z^2\rightharpoonup z=\frac{1}{2}(1+g+h-gh)\rightharpoonup z=\frac{1}{2}(z+ghz+gz-hz)$$ and Eq.~\eqref{xishu}, we simultaneously obtain that
$a_5a_6+a_6a_8+a_7a_5+a_8a_7=\frac{1}{2}$ and $a_5a_7+a_6a_5+a_7a_8+a_8a_6=-\frac{1}{2}$, which is a contradiction.
\end{remark}

Now for each of the left actions given in TABLEs \ref{H8}--\ref{H8--}, it is straightforward to check that the linear map $\leftharpoonup$ given by Eq.~\eqref{eq:r-action} is a right module coalgebra action on $H_8$.
Then according to Theorem~\ref{thm:construct-mp},  we obtain the following  classification theorem.
\begin{theorem}\label{thm:mpa_H_8}
There are 6 matched pairs of actions on $H_8$, and they are defined by the left actions $\rightharpoonup$ given in TABLEs \ref{H8}--\ref{H8--} respectively.
\end{theorem}

\section{Matched pairs of actions on $H_8$ irrelevant to its coquasitriangular structures}\label{sec:ir1}

Recall that a {\bf coquasitriangular Hopf algebra} is a Hopf algebra $H$ with a convolution-invertible bilinear map $\calr:H\otimes H\to\bk$ such that
\begin{eqnarray*}
\calr(a_1\otimes b_1)a_2b_2&=&b_1a_1\calr(a_2\otimes b_2),\\
\calr(a\otimes bc)&=&\calr(a_1\otimes c)\calr(a_2\otimes b),\\
\calr(ab\otimes c)&=&\calr(a\otimes c_1)\calr(b\otimes c_2)
\end{eqnarray*}
for any $a,b,c\in H$~\cite[\S~2.2]{Ma}.

In \cite{FS}, the authors pointed out that any coquasitriangular Hopf algebra $(H,\calr)$ naturally provides a matched pair of actions $(H,\rightharpoonup,\leftharpoonup)$ defined by
\begin{eqnarray}\label{laction}
a\rightharpoonup b&=& \calr^{-1}(a_1\otimes b_1)b_2\calr(a_2\otimes b_3),
\end{eqnarray}
\begin{eqnarray}\label{raction}
a\leftharpoonup b&=& \calr^{-1}(a_1\otimes b_1)a_2\calr(a_3\otimes b_2).
\end{eqnarray}

So far we have classified all matched pairs of actions on $H_8$. In this section, we pick out those not derived from coquasitriangular structures on $H_8$. In order to do so, we need to review the following results in a more general context.

Suzuki's Hopf algebras $A_{NL}^{v\lambda}$ were introduced in \cite{Su} as a family of finite dimensional cosemisimple Hopf algebras generated by a comatrix basis of the $2\times 2$-matrices. In particular, the algebra $A_{12}^{+-}$ in this family is exactly the Kac-Paljutkin Hopf algebra $H_8$ with the following different presentation. It is generated by 4 elements $x_{11}$, $x_{12}$, $x_{21}$ and $x_{22}$ with the relations
$$x_{11}^2=x_{22}^2,\quad x_{12}^2=x_{21}^2,\quad x_{11}^2+x_{12}^2=1,\quad x_{21}x_{12}=-x_{12}x_{21},\quad x_{11}x_{22}=x_{22}x_{11},$$
$$x_{11}x_{12}=x_{12}x_{11}=x_{11}x_{21}=x_{21}x_{11}=x_{12}x_{22}=x_{22}x_{12}=x_{21}x_{22}=x_{22}x_{21}=0.$$
A linear basis of $A_{12}^{+-}$ is given by $\{1,x_{11},x_{12},x_{21},x_{22},x_{11}^2,x_{11}x_{22},x_{12}x_{21}\}$. The coalgebra structure and the antipode are defined by
$$\Delta(x_{ij})=x_{i1}\otimes x_{1j}+x_{i2}\otimes x_{2j},\quad \varepsilon(x_{ij})=\delta_{ij},\quad S(x_{ij})=x_{ji},\quad 1\leq i,j\leq 2.$$
The set of group-like elements of $A_{12}^{+-}$ is
$$G(A_{12}^{+-})=\{1, x_{11}^2-x_{12}^2, x_{11}x_{22}+\sqrt{-1}x_{12}x_{21}, x_{11}x_{22}-\sqrt{-1}x_{12}x_{21}\}.$$

Note that all coquasitriangular structures of the Suzuki Hopf algebra $A_{NL}^{v\lambda}$ have been determined in \cite{Su} (see also~\cite[Theorem~5.7]{Wa}).
In particular, there are two families of coquasitriangular structures on the Kac-Paljutkin algebra $A_{12}^{+-}$ given by
$$\{\sigma_{\alpha,\,\beta}\,|\,\alpha,\,\beta\in\bk^{\times},\alpha\beta=1,(\alpha\beta^{-1})^2=-1\}\cup
\{\tau_{\gamma,\,\xi}\,|\,\gamma,\,\xi\in\bk^{\times}, \gamma^2=\xi^2=1\},$$
where
\begin{equation}\label{coquasitriangular1}
\sigma_{\alpha,\,\beta}
\begin{pmatrix}
x_{11}\otimes x_{11}& x_{11}\otimes x_{12}& x_{11}\otimes x_{21}& x_{11}\otimes x_{22}\\
x_{12}\otimes x_{11}& x_{12}\otimes x_{12}& x_{12}\otimes x_{21}& x_{12}\otimes x_{22}\\
x_{21}\otimes x_{11}& x_{21}\otimes x_{12}& x_{21}\otimes x_{21}& x_{21}\otimes x_{22}\\
x_{22}\otimes x_{11}& x_{22}\otimes x_{12}& x_{22}\otimes x_{21}& x_{22}\otimes x_{22}
\end{pmatrix}=\begin{pmatrix}
0& 0& 0& 0\\
0& \alpha& \beta & 0\\
0& \beta& \alpha& 0\\
0& 0& 0& 0
\end{pmatrix},
\end{equation}
\begin{equation}\label{coquasitriangular2}
\tau_{\gamma,\,\xi}
\begin{pmatrix}
x_{11}\otimes x_{11}& x_{11}\otimes x_{12}& x_{11}\otimes x_{21}& x_{11}\otimes x_{22}\\
x_{12}\otimes x_{11}& x_{12}\otimes x_{12}& x_{12}\otimes x_{21}& x_{12}\otimes x_{22}\\
x_{21}\otimes x_{11}& x_{21}\otimes x_{12}& x_{21}\otimes x_{21}& x_{21}\otimes x_{22}\\
x_{22}\otimes x_{11}& x_{22}\otimes x_{12}& x_{22}\otimes x_{21}& x_{22}\otimes x_{22}
 \end{pmatrix}=\begin{pmatrix}
\gamma& 0& 0& \xi\\
0& 0& 0& 0\\
0& 0& 0& 0\\
-\xi& 0& 0& \gamma
\end{pmatrix},
\end{equation}
whose convolution inverse are $\sigma^{-1}_{\alpha,\,\beta}=\sigma_{\alpha^{-1},\,\beta^{-1}}$ and $\tau^{-1}_{\gamma,\, \xi}=\tau_{\gamma^{-1},\,\xi^{-1}}$.
Then we can get four left actions through Eq.~\eqref{laction}, among which two
in TABLEs \ref{A12} and \ref{A12'} are obtained from the first family of coquasitriangular structures \eqref{coquasitriangular1} and the other two in TABLEs \ref{A12''} and \ref{A12'''} are obtained from the second family of coquasitriangular structures \eqref{coquasitriangular2}.
\begin{table}[H]
\centering
\begin{minipage}{0.46\textwidth}
\centering
\caption{The 1st left action on $A_{12}^{+-}$ }\label{A12}
\begin{tabular}{c|cccc}
$\rightharpoonup$&$x_{11}$&$x_{12}$&$x_{21}$&$x_{22}$\\
\hline
$x_{11}$&$x_{22}$&$\sqrt{-1}x_{21}$&$-\sqrt{-1}x_{12}$&$x_{11}$\\
$x_{12}$&0&0&0&0\\
$x_{21}$&0&0&0&0\\
$x_{22}$&$x_{22}$&$-\sqrt{-1}x_{21}$&$\sqrt{-1}x_{12}$&$x_{11}$\\
\end{tabular}
\end{minipage}
\hspace{0.05\textwidth} 
\begin{minipage}{0.46\textwidth}
\centering
\caption{The 2nd left action on $A_{12}^{+-}$}\label{A12'}
\begin{tabular}{c|cccc}
$\rightharpoonup$&$x_{11}$&$x_{12}$&$x_{21}$&$x_{22}$\\
\hline
$x_{11}$&$x_{22}$&$-\sqrt{-1}x_{21}$&$\sqrt{-1}x_{12}$&$x_{11}$\\
$x_{12}$&0&0&0&0\\
$x_{21}$&0&0&0&0\\
$x_{22}$&$x_{22}$&$\sqrt{-1}x_{21}$&$-\sqrt{-1}x_{12}$&$x_{11}$\\
\end{tabular}
\end{minipage}
\end{table}

\begin{table}[H]
\centering
\begin{minipage}{0.47\textwidth}
\centering
\caption{The 3rd left action on $A_{12}^{+-}$}\label{A12''}
\begin{tabular}{c|cccc}
$\rightharpoonup$&$x_{11}$&$x_{12}$&$x_{21}$&$x_{22}$\\
\hline
$x_{11}$&$x_{11}$&$-x_{12}$&$-x_{21}$&$x_{22}$\\
$x_{12}$&0&0&0&0\\
$x_{21}$&0&0&0&0\\
$x_{22}$&$x_{11}$&$x_{12}$&$x_{21}$&$x_{22}$\\
\end{tabular}
\end{minipage}
\hspace{0.04\textwidth} 
\begin{minipage}{0.47\textwidth}
\centering
\caption{The 4th left action on $A_{12}^{+-}$}\label{A12'''}
\begin{tabular}{c|cccc}
$\rightharpoonup$&$x_{11}$&$x_{12}$&$x_{21}$&$x_{22}$\\
\hline
$x_{11}$&$x_{11}$&$x_{12}$&$x_{21}$&$x_{22}$\\
$x_{12}$&0&0&0&0\\
$x_{21}$&0&0&0&0\\
$x_{22}$&$x_{11}$&$-x_{12}$&$-x_{21}$&$x_{22}$\\
\end{tabular}
\end{minipage}
\end{table}
Next we identify the aforementioned two presentations of $H_8$ whose sets of generators are $\{g,h,z\}$ and $\{x_{ij}\}_{i,j=1,2}$ respectively. First we take
$$\quad g=x_{11}x_{22}+\sqrt{-1}x_{12}x_{21}, \quad h=x_{11}x_{22}-\sqrt{-1}x_{12}x_{21},$$
and then give an expansion $z=a_1x_{11}+a_2x_{12}+a_3x_{21}+a_4x_{22}$.
On the one hand,
$$z^2=(a_1x_{11}+a_2x_{12}+a_3x_{21}+a_4x_{22})^2=(a_1^2+a_4^2)x_{11}^2+(a_2^2+a_3^2)x_{12}^2+2a_1a_4x_{11}x_{22}.$$
On the other hand,
$$z^2=\frac{1}{2}(1+g+h-gh)=x_{12}^2+x_{11}x_{22}.$$
By the comparison of coefficients, we have
\begin{equation}\label{fangchengzu4}
\begin{cases}
a_1^2+a_4^2=0,\\
a_2^2+a_3^2=1,\\
a_1a_4=\frac{1}{2}.
\end{cases}
\end{equation}
Moreover, since
\begin{align*}
\Delta(z)&=\frac{1}{2} (z \otimes z + gz \otimes z+ z \otimes hz- gz \otimes hz )\\
&=\Delta(a_1x_{11}+a_2x_{12}+a_3x_{21}+a_4x_{22})\\
&=a_1(x_{11}\otimes x_{11}+x_{12}\otimes x_{21})+a_2(x_{11}\otimes x_{12}+x_{12}\otimes x_{22})\\
&\quad+a_3(x_{21}\otimes x_{11}+x_{22}\otimes x_{21})+a_4(x_{21}\otimes x_{12}+x_{22}\otimes x_{22}),
\end{align*}
we obtain 16 equations, among which the following are the crucial ones.
\begin{equation}\label{fangchengzu5}
\begin{cases}
a_1^2+a_4a_1+a_1a_4-a_4^2=2a_1,\\
a_1a_2+a_4a_2-a_1a_3\sqrt{-1}+a_4a_3\sqrt{-1}=2a_2,\\
a_2a_3+a_3^2\sqrt{-1}+a_2^2\sqrt{-1}+a_3a_2=2a_1,\\
a_4^2+a_1a_4+a_4a_1-a_1^2=2a_4.
\end{cases}
\end{equation}
Combining Eqs.~\eqref{fangchengzu4} and \eqref{fangchengzu5}, we confirm that $$z=\frac{1+\sqrt{-1}}{2}x_{11}+\frac{1-\sqrt{-1}}{2}x_{22}+\frac{\sqrt{2}}{2}x_{12}+\frac{\sqrt{2}}{2}x_{21}.$$

To relate TABLEs \ref{A12}--\ref{A12'''}
with 4 ones among TABLEs \ref{H8}--\ref{H8--} for the left actions of $H_8$ on itself, we calculate some key elements within them. Taking TABLE \ref{A12} as an example, we have
\begin{align*}
g\rightharpoonup z&=(x_{11}x_{22}+\sqrt{-1}x_{12}x_{21})\rightharpoonup(\frac{1+\sqrt{-1}}{2}x_{11}+\frac{1-\sqrt{-1}}{2}x_{22}+\frac{\sqrt{2}}{2}x_{12}+\frac{\sqrt{2}}{2}x_{21})\\
&=\frac{1+\sqrt{-1}}{2}x_{11}+\frac{1-\sqrt{-1}}{2}x_{22}-\frac{\sqrt{2}}{2}x_{12}-\frac{\sqrt{2}}{2}x_{21}\\
&=ghz.\\[.5em]
h\rightharpoonup z&=(x_{11}x_{22}-\sqrt{-1}x_{12}x_{21})\rightharpoonup(\frac{1+\sqrt{-1}}{2}x_{11}+\frac{1-\sqrt{-1}}{2}x_{22}+\frac{\sqrt{2}}{2}x_{12}+\frac{\sqrt{2}}{2}x_{21})\\
&=\frac{1+\sqrt{-1}}{2}x_{11}+\frac{1-\sqrt{-1}}{2}x_{22}-\frac{\sqrt{2}}{2}x_{12}-\frac{\sqrt{2}}{2}x_{21}\\
&=ghz.\\[.5em]
z\rightharpoonup z&=(\frac{1+\sqrt{-1}}{2}x_{11}+\frac{1-\sqrt{-1}}{2}x_{22}+\frac{\sqrt{2}}{2}x_{12}+\frac{\sqrt{2}}{2}x_{21})\\
&\quad\rightharpoonup(\frac{1+\sqrt{-1}}{2}x_{11}+\frac{1-\sqrt{-1}}{2}x_{22}+\frac{\sqrt{2}}{2}x_{12}+\frac{\sqrt{2}}{2}x_{21})\\
&=\frac{1-\sqrt{-1}}{2}x_{11}+\frac{1+\sqrt{-1}}{2}x_{22}+\frac{\sqrt{2}}{2}x_{12}-\frac{\sqrt{2}}{2}x_{21}\\
&=\frac{1-\sqrt{-1}}{2}gz+\frac{1+\sqrt{-1}}{2}hz.
\end{align*}

Therefore, we observe that the left action in TABLE \ref{A12} is consistent with that in TABLE \ref{H8}.
Similarly, we check that the left actions in TABLEs \ref{A12'}--\ref{A12'''} are successively consistent with those in TABLEs \ref{H8'}--\ref{H8'''}. In other words, we obtain the following result.
\begin{theorem}\label{thm:mpa_H_8_coqt}
There are two matched pairs of actions on $H_8$ not derived from its coquasitriangular structures, and they are defined by the left actions $\rightharpoonup$
given in TABLEs \ref{H8-} and \ref{H8--} respectively.
\end{theorem}


\section{The Yang-Baxter operators associated to matched pairs of actions on $H_8$}\label{sec:ir2}

Let $(H,\rightharpoonup,\leftharpoonup)$ be a matched pair of actions on a Hopf algebra $H$. According to~\cite[Theorem~5.27]{GGV} or \cite[Theorem~3.2]{Li},
the linear operator $r:H\otimes H\to H\otimes H$ defined by
\begin{eqnarray}\label{eq:pyb}
r(x\otimes y) &=& (x_1\rightharpoonup y_1)\otimes (x_2\leftharpoonup y_2),\quad\forall x,y\in H
\end{eqnarray}
is a bijective solution of the braid equation:
$$(r \otimes \id)(\id \otimes r)(r \otimes \id) = (\id \otimes r)(r \otimes \id)(\id \otimes r).$$
Namely, $r$ is a Yang-Baxter operator.

In this section, we study the Yang-Baxter operators associated to those matched pairs of actions on $H_8$ classified in Section~\ref{sec:H8}. First we reorganize the results given in \cite[Theorem~3.5, Corollary~3.7]{Li} as follows.

\begin{theorem}[\cite{Li}]\label{th:duihe}
Let $(H,\rightharpoonup,\leftharpoonup)$ be a matched pair of actions on a Hopf algebra H. The Yang-Baxter operator $r$ defined by Eq.~\eqref{eq:pyb} is involutive, i.e. $r^2=\id^{\otimes 2}$, if and only if
\begin{equation}\label{duihe1}
x\leftharpoonup y=S(x_1\rightharpoonup y)\rightharpoonup x_2,\quad \forall x,y\in H.
\end{equation}
In this situation, we also have
\begin{equation}\label{duihe2}
S(x\leftharpoonup y)=S(y)\rightharpoonup S(x),\quad \forall x,y\in H.
\end{equation}
\end{theorem}

By calculation, we find that the 4 matched pairs of actions on $H_8$ defined by the left actions $\rightharpoonup$ in TABLEs \ref{H8}--\ref{H8'''} respectively do not satisfy Eq.~\eqref{duihe2}. Take the one from TABLE \ref{H8} as an example.
On the one hand,
\begin{align*}
z\leftharpoonup gz&\stackrel{\eqref{eq:r-action}}=\frac{1}{4}(S(z\rightharpoonup gz)zgz+S(z\rightharpoonup z)zgz+S(z\rightharpoonup gz)zghz-S(z\rightharpoonup z)zghz\\
&\quad+S(gz\rightharpoonup gz)zgz+S(gz\rightharpoonup z)zgz+S(gz\rightharpoonup gz)zghz-S(gz\rightharpoonup z)zghz\\
&\quad+S(z\rightharpoonup gz)hzgz+S(z\rightharpoonup z)hzgz+S(z\rightharpoonup gz)hzghz-S(z\rightharpoonup z)hzghz\\
&\quad-S(gz\rightharpoonup gz)hzgz-S(gz\rightharpoonup z)hzgz-S(gz\rightharpoonup gz)hzghz+S(gz\rightharpoonup z)hzghz)\\
&=\frac{1+\sqrt{-1}}{2}gz+\frac{1-\sqrt{-1}}{2}hz,
\end{align*}
so we have $$S(z\leftharpoonup gz)=S(\frac{1+\sqrt{-1}}{2}gz+\frac{1-\sqrt{-1}}{2}hz)=\frac{1+\sqrt{-1}}{2}hz+\frac{1-\sqrt{-1}}{2}gz.$$
On the other hand,
$$S(gz)\rightharpoonup S(z)=zg\rightharpoonup z=\frac{1-\sqrt{-1}}{2}hz+\frac{1+\sqrt{-1}}{2}gz.$$
That is, Eq.~\eqref{duihe2} does not hold.

By contrast, for the 2 matched pairs of actions on $H_8$ defined by the left actions $\rightharpoonup$ in TABLEs \ref{H8-} and \ref{H8--} respectively, we can verify that they satisfy Eq.~\eqref{duihe1}. For the sake of completeness, we list their corresponding right actions $\leftharpoonup$ as follows.

\begin{table}[H]
\centering
\begin{minipage}{0.45\textwidth}
\centering
\caption{The right action \hspace*{4em} on $H_8$ in Case 2.1}
\renewcommand{\arraystretch}{1.2}
{\footnotesize\begin{tabular}{|c|c|c|c|c|c|c|c|c|}
\hline
$\leftharpoonup$&1&$g$&$h$&$gh$&$z$&$gz$&$hz$&$ghz$\\
\hline
1&1&$1$&$1$&$1$&$1$&$1$&$1$&$1$\\
\hline
$g$&$g$&$g$&$g$&$g$&$h$&$h$&$h$&$h$\\
\hline
$h$&$h$&$h$&$h$&$h$&$g$&$g$&$g$&$g$\\
\hline
$gh$&$gh$&$gh$&$gh$&$gh$&$gh$&$gh$&$gh$&$gh$\\
\hline
$z$&$z$&$z$&$z$&$z$&$z^{3}$&$z^{3}$&$z^{3}$&$z^{3}$\\
\hline
$gz$&$gz$&$gz$&$gz$&$gz$&$hz^{3}$&$hz^{3}$&$hz^{3}$&$hz^{3}$\\
\hline
$hz$&$hz$&$hz$&$hz$&$hz$&$gz^{3}$&$gz^{3}$&$gz^{3}$&$gz^{3}$\\
\hline
$ghz$&$ghz$&$ghz$&$ghz$&$ghz$&$ghz^{3}$&$ghz^{3}$&$ghz^{3}$&$ghz^{3}$\\
\hline
\end{tabular}}
\end{minipage}
\hspace{0.08\textwidth}
\begin{minipage}{0.45\textwidth}
\centering
\caption{The right action \hspace*{4em} on $H_8$ in Case 2.2}
\renewcommand{\arraystretch}{1.2}
{\footnotesize\begin{tabular}{|c|c|c|c|c|c|c|c|c|}
\hline
$\leftharpoonup$&1&$g$&$h$&$gh$&$z$&$gz$&$hz$&$ghz$\\
\hline
1&1&$1$&$1$&$1$&$1$&$1$&$1$&$1$\\
\hline
$g$&$g$&$g$&$g$&$g$&$h$&$h$&$h$&$h$\\
\hline
$h$&$h$&$h$&$h$&$h$&$g$&$g$&$g$&$g$\\
\hline
$gh$&$gh$&$gh$&$gh$&$gh$&$gh$&$gh$&$gh$&$gh$\\
\hline
$z$&$z$&$z$&$z$&$z$&$ghz^{3}$&$ghz^{3}$&$ghz^{3}$&$ghz^{3}$\\
\hline
$gz$&$gz$&$gz$&$gz$&$gz$&$gz^{3}$&$gz^{3}$&$gz^{3}$&$gz^{3}$\\
\hline
$hz$&$hz$&$hz$&$hz$&$hz$&$hz^{3}$&$hz^{3}$&$hz^{3}$&$hz^{3}$\\
\hline
$ghz$&$ghz$&$ghz$&$ghz$&$ghz$&$z^{3}$&$z^{3}$&$z^{3}$&$z^{3}$\\
\hline
\end{tabular}}
\end{minipage}
\end{table}

So according to Theorem \ref{th:duihe}, we have the following result.
\begin{theorem}\label{thm:mpa_H_8_ybo}
The Yang-Baxter operators associated to matched pairs of actions on $H_8$ which are defined by the left actions $\rightharpoonup$ given in TABLEs \ref{H8-} and \ref{H8--} respectively are involutive.
\end{theorem}

\vspace{0.1cm}
 \noindent
{\bf Acknowledgements.} This work is supported by National Natural Science Foundation of China (12071094, 12171155), and Basic and Applied Basic Research Foundation of Guangdong Province (2022A1515010357).

\bibliographystyle{amsplain}

\begin{thebibliography}{99}


\bibitem{Ag} A. L. Agore, Classifying bicrossed products of two Taft algebras, \textit{J. Pure Appl. Algebra} {\bf 222} (2018), 914--930.

\bibitem{ABM} A. L. Agore, C. G. Bontea and G. Militaru, Classifying bicrossed products of Hopf algebras, \textit{Algebr. Represent. Theory} {\bf 17} (2014), 227--264.

\bibitem{AGV} I. Angiono, C. Galindo and L. Vendramin, Hopf braces and Yang-Baxter operators, \textit{Proc. Amer. Math. Soc.} {\bf 145} (2017), 1981--1995.

\bibitem{Bo} C. G. Bontea, Classifying bicrossed products of two Sweedler's Hopf algebras, \textit{Czech. Math. J.} {\bf 64} (2014), 419--431.


\bibitem{FS}
D. Ferri and A. Sciandra, Matched pairs and Yetter-Drinfeld braces, arXiv:2406.10009.

\bibitem{GGV}
J. A. Guccione, J. J. Guccione and C. Valqui, Set-theoretic type solutions of the braid equation, arXiv:2008.13494 (v5); with the original version already published in \textit{J. Algebra}


\bibitem{HZ1} J. Hu and Y. Zhang, The $\beta$-character algebra and a commuting pair in Hopf algebras, \textit{Algebr. Represent. Theory} \textbf{10} (2007), 497--516.

\bibitem{KP} G. I. Kac and V. G. Paljutkin, Finite ring groups, {\it Trudy Moskov. Mat. Ob$\check{s}\check{c}$} {\bf 15} (1966), 224--261.

\bibitem{Li} Y. Li, Matched pairs and Yang-Baxter operators, arXiv:2501.11975.


\bibitem{LNW} D. Lu, Y. Ning and D. Wang, The bicrossed products of $H_4$ and $H_8$, \textit{Czech. Math. J.} {\bf 70} (2020), 959--977.

\bibitem{LYZ} J. Lu, M. Yan and Y. Zhu,  On the set-theoretical Yang-Baxter equation, \textit{Duke Math. J.} {\bf 104} (2000),   1--18.

\bibitem{Ma2} S. Majid, Physics for algebraists: Non-commutative and non-cocommutative Hopf algebras by a bicrossproduct construction, \textit{J. Algebra} {\bf 130} (1990), 17--64.

\bibitem{Ma} S. Majid, Foundations of quantum group theory, Cambridge University Press, 1995.

\bibitem{Mas}
A.~Masuoka, Semisimple Hopf algebras of dimension {$6,8$}, \textit{Israel J. Math.} \textbf{92} (1995), 361--373.

\bibitem{Mon} S. Montgomery, Hopf algebras and Their Actions on Rings, Amer. Math. Soc., Regional Conf. Ser. in Math., \textbf{82}, 1993.

\bibitem{Pa} D. Pansera, A class of semisimple Hopf algebras acting on quantum polynomial algebras, in: Rings, modules and codes, 303--316, Contemp. Math., {\bf 727.}, 2019.


\bibitem{SV} D.~S.~Sage and M.~D.~Vega, Twisted Frobenius-Schur indicators for Hopf algebras, \textit{J. Algebra} \textbf{354} (2012), 136--147.


\bibitem{Shi} Y. Shi, Finite dimensional Hopf algebras over the Kac-Paljutkin algebra $H_8$, \textit{Rev. Un. Mat. Argentina} \textbf{60} (2019), 265--298.


\bibitem{Si} W. M. Singer, Extension theory for connected Hopf algebras.  \textit{J. Algebra} \textbf{21} (1972), 1--16.

\bibitem{St} D. Stefan, Hopf algebras of low dimension, \textit{J. Algebra} {\bf 211} (1999), 343--361.

\bibitem{Su} S. Suzuki, A family of braided cosemisimple Hopf algebras of finite dimension, \textit{Tsukuba J. Math.} {\bf 22} (1998), 1--29.

\bibitem{Ta} M. Takeuchi, Matched pairs of groups and bismash products of Hopf algebras, \textit{Comm. Algebra}  {\bf 9} (1981), 841--882.


\bibitem{Wa} M. Wakui, Polynomial invariants for a semisimple and cosemisimple Hopf algebra of finite dimension, \textit{J. Pure Appl. Algebra} {\bf 214} (2010), 701--728

\bibitem{ZL} K. Zhou, G. Liu, On the quasitriangular structures of abelian extensions of $\mathbb Z_2$, \textit{Comm. Algebra} {\bf 49} (2021), 4755--4762.

\end{thebibliography}

\end{document}